\documentclass{article}

\usepackage{longtable}
\usepackage{setspace}

\begin{document}
\title{A Search for Spreads of Hermitian Unitals}
\author{Jeremy M. Dover\\
Dover Networks LLC}
\date{February 3, 2006, updated February 4, 2017}
\maketitle

\begin{abstract}
A {\em spread} of a Hermitian unital in $PG(2,q^2)$ is a set of $q^2-q+1$ 
pairwise disjoint blocks that partition the points of the unital. In this 
paper, we discuss the results of an exhaustive computer search for spreads 
of Hermitian unitals of small orders, namely the Hermitian unitals in 
$PG(2,16)$, $PG(2,25)$ and $PG(2,49)$. 
\end{abstract}

\section{Introduction}
Let $q$ be a prime power, and let $\Pi$ be the projective plane over $GF(q)$. 
A {\it Hermitian unital} is the set of absolute points of a unitary polarity 
in $\Pi$. A Hermitian unital $U$ has the property that any line of $\Pi$ meets 
$U$ in either 1 or $q+1$ points, and these lines are called {\it tangents} 
and {\it secants}, respectively. Considering the intersections of $U$ and its 
secants as blocks, the Hermitian unital in $PG(2,q^2)$ is a $2-(q^3+1,q+1,1)$ 
design. When no confusion can arise, we will refer to the set of lines of 
$\Pi$ whose intersections with $U$ form a spread of $U$ as a spread as well.

An interesting question with any design is whether or not one can partition 
the pointset of the design into pairwise disjoint blocks, creating a 
{\it spread}. Spreads of many different types of spaces have been studied, and 
in~\cite{dover:ree}, the author investigated spreads of Ree unitals.

The obvious method for creating a spread of a Hermitian unital $U$ in $\Pi$ 
requires one to start with a point $P \in \Pi$ outside of $U$. The set of 
secants to $U$ through $P$, together with $P^\perp$, the polar of $P$ under 
the unitary polarity associated with $U$, form a spread of $U$. We refer to 
such a spread as being {\it regular}.

In addition to the regular spread, Baker, et. al.~\cite{beks} provide several 
examples of spreads of Hermitian unitals. The focus in that work is on 
finding spreads which are {\it orthogonally divergent}, i.e., spreads such 
that no line of the spread contains the pole of any other line of the spread, 
with respect to the polarity associated with the Hermitian unital.

The first type of spread from~\cite{beks} is a cyclic spread, i.e., a spread 
which can be obtained by taking the orbit of one of its members under a cyclic 
group of order $q^2-q+1$. Baker, et. al., prove that a cyclic spread of a 
Hermitian unital in $PG(2,q^2)$ exists if and only if $q$ is even.

The second type of spread given in~\cite{beks} involves an interesting 
partition of the Hermitian unital. Let $U$ be the Hermitian unital defined by 
the equation $x^{q+1}+y^{q+1}+z^{q+1} = 0$, where $x$, $y$, and $z$ are 
homogeneous coordinates for $\Pi=PG(2,q^2)$. Let $P=(1,0,0)$, $Q=(0,1,0)$ and 
$R=(0,0,1)$ be the vertices of a self-polar triangle with respect to $U$. 
Then $U$ is partitioned by the sides of the triangle $\triangle PQR$ and 
the sets

\begin{displaymath}
T_a = \{(1,y,z):y^{q+1} = -(1+a), z^{q+1} = a\},
\end{displaymath}
for each $a \in GF(q) \setminus \{0,-1\}$. Each of the $T_a$'s is triply-ruled,
in that there are sets of $q+1$ secants to $U$ through each of $P$, $Q$ and 
$R$ whose intersection with $U$ exactly cover $T_a$. These sets are:
\begin{eqnarray}
H_a & = & \{[u,0,1]:u^{q+1} = a\}\nonumber\\
V_a & = & \{[v,1,0]:v^{q+1} = -(1+a)\}\label{tripleruling}\\
D_a & = & \{[0,w,1]:w^{q+1} = -a/(1+a)\nonumber\}
\end{eqnarray}

Hence one can create a spread of $U$ by selecting one of the three ruling 
families $H_a$, $V_a$, or $D_a$ for each $T_a$ and adjoining the three sides 
of the triangle $\triangle PQR$.

Based on the data generated in the search described here, the
author~\cite{dover:subreg} provides a generalization
of this construction to create a family of {\em subregular} spreads of
the Hermitian unital. Using an analogy with the spreads of $PG(3,q)$,
we call the spreads generated in Baker, et. al.~\cite{beks} an {\em Andr\'e spread},
noting that there exist subregular spreads of the Hermitian unital which are
not Andr\'e spreads.

\section{Search Procedure}

Our basic algorithm for searching for spreads is a standard backtrack, where we
sequentially attempt to add lines to a partial spread until we either get to a 
full spread or have no lines remaining that are disjoint from our partial
spread, in which case we backtrack. Despite multiple attempts to increase the 
efficiency of this backtrack, we felt it unlikely that it would be able to 
complete in a reasonable amount of time in $PG(2,49)$. Moreover, this backtrack
produces an enormous number of isomorphic spreads which must then be sorted.

To improve the performance of this procedure, we perform an aggressive 
isomorph rejection phase prior to backtracking to determine a number of
``starter'' configurations that exhaustively cover the search space as follows.
We set up our initial data structures using Magma~\cite{cannon}, including the 
unital, defined to be $U=\{(x,y,z):x^{q+1}+y^{q+1}+z^{q+1} = 0\}$, and we fix
the secant $[1,0,0]$. Since the automorphism
group of the Hermitian unital is transitive on secants, we may assume that
$[1,0,0]$ is in our spread.

We define a {\it starter configuration} $S$ to be a pair ($S_i$,$S_e$) of sets
such that $S_i$ is the set of secants to the unital which are assumed to be 
in the spread we will be constructing, and $S_e$ is the set of secants 
specifically excluded from this potential spread. Thus our initial
starter configuration is $(\{[1,0,0]\},\emptyset)$.

From this initial starter configuration $S=(\{[1,0,0]\},\emptyset)$, we may 
generate a set of new starter 
configurations which also cover the search space as follows. Compute the 
stabilizer $L$ within the automorphism group of the unital of the set 
$\{[1,0,0]\}$. Let $\{O_j:1 \leq j \leq n\}$ be the orbits under $L$ of those 
secants which do not meet $[1,0,0]$. 
Letting $\ell_j$ be an orbit representative of $O_j$, it
is clear that every spread of $U$ contains a set of lines isomorphic to
$\{[1,0,0],\ell_j\}$ for some $1 \leq j \leq n$. Thus we can partition our search space into $n$ new
search spaces, namely $(\{[1,0,0],\ell_j\},\emptyset)$ for each $1 \leq j \leq n$.

Partitioning the search space in this way does not give us maximum
advantage. Suppose there exists a spread ${\cal S}$ of $U$ that contains the
lines $\{[1,0,0],m_j,m_k\}$ where $m_j \in O_j$ and $m_k \in O_k$ with 
$j \neq k$. Then our search will find spreads isomorphic to ${\cal S}$ in 
each of the search spaces with starters 
$(\{[1,0,0],\ell_j\},\emptyset)$ and $(\{[1,0,0],\ell_k\},\emptyset)$. To be 
more efficient, we could have the search space with starter 
$(\{[1,0,0],\ell_j\},\emptyset)$ exclude the blocks from the orbit $O_k$ and 
still find a spread isomorphic to ${\cal S}$. Therefore, we will instead 
partition our search space $(\{[1,0,0]\},\emptyset)$ into the search spaces 
with starter configurations:
\begin{displaymath}
\left( \{[1,0,0],\ell_j\},\bigcup_{k = j+1}^{n} O_k\right)
\end{displaymath}

We can iterate this procedure with each of the new starters
to further decompose our search space. Let $(S_i,S_e)=(\{[1,0,0],\ell_j\},
\bigcup_{k=j+1}^n O_k)$ be one of our new starter configurations. Again 
compute the setwise stabilizer $L'$ of $S_i$ in the unital $U$. We look 
at the orbits $O'_k$, $1 \leq k \leq m$ of secants under $L'$. As before we
need not consider those orbits whose secants meet a secant of $S_i$ in
a point of $U$. However, we may also exclude from consideration all of those
orbits which contain a secant in $S_e$, even if those orbits contain secants
which are not in $S_e$.

To see why this exclusion is valid, suppose there exists a spread ${\cal S}'$
which contains the secants $[1,0,0]$, $\ell_j$, and $\ell'$, where $\ell' 
\notin S_e$, but there exists a line $\ell'' \in S_e$ in $\ell'$'s orbit under
$L'$. Then there exists an automorphism of $U$ in $L'$ that maps $\{[1,0,0],
\ell_j\}$ onto itself and maps $\ell'$ onto $\ell''$, and thus the spread
${\cal S}'$ onto a new spread ${\cal S}''$ containing $\{[1,0,0], \ell_j,
\ell''\}$. Since $\ell''$ is in $S_e$, it must be in an orbit $O_m$ with
$m > j$, which means that a spread isomorphic to ${\cal S}''$, and hence 
${\cal S}'$ will be found in the search space with starter configuration 
$(\{[1,0,0],\ell_m\},\bigcup_{k=m+1}^n O_k)$.

We use Magma to iterate this procedure as long as one of the starter 
configurations has an initial partial spread with a nontrivial stabilizer in
the automorphism group of $U$, or the number of starter configurations becomes
too unwieldy. At the end of the isomorph rejection phase, the
searcher will be left with a number of starter configurations which are the
seeds to our backtracking algorithm.

Using Magma, we set up the data structures for our C backtracking 
algorithm. For each starter configuration $(S_i,S_e)$, we create the initial 
partial spread from $S_i$ and a list of all candidate secants that can 
be used to complete the initial partial spread to a full spread. This 
set consists of those secants which neither meet an element of $S_i$ 
in a point of the unital nor are in the set of excluded secants, $S_e$.

After all of the backtrack searches are completed, we use a number of Perl
scripts to marshall the spreads together, leveraging Magma to compute several
combinatorial invariants as an initial isomorphism sort. We then pass 
completely back into Magma to execute a full isomorphism check and produce
a comprehensive list of projectively inequivalent spreads.

Using a custom-built Pentium 4 3.0 GHz machine with 2GB of RAM running RedHat 
Linux 9.0, the isomorph rejection phase of the search ran almost instantaneous 
in PG(2,16), took a few minutes for PG(2,25), and lasted approximately 20 hours for 
PG(2,49). The backtracking phase was again almost instantaneous for PG(2,16),
less than a minute for PG(2,25), and lasted approximately 10 weeks for
PG(2,49). In the PG(2,49) case, we searched a total of 20,846 starter 
configurations, many trivial due to large numbers of excluded secants. As none 
of the individual backtrack runs took more than 26 hours to complete, the
highly parallelizable nature of our search suggests that an exhaustive search 
for spreads of the Hermitian unital in PG(2,64) is not beyond reach
with more robust hardware.

\section{$PG(2,16)$}

Let $w$ be a primitive element of $GF(16)$ with minimal polynomial $x^4+x+1$
over $GF(2)$. Our search found exactly three spreads of the Hermitian unital
$U = \{(x,y,z):x^{q+1}+y^{q+1}+z^{q+1} = 0\}$ in $PG(2,16)$, all of which were 
previously known, namely the regular spread, as well as the cyclic 
spread and a spread from the triply-ruled partition as constructed in Baker, 
et. al.~\cite{beks}.

Constructions and properties of these spreads are presented in 
Table~\ref{PG216table}. In order to report the lines of each spread,
we start with the regular spread $S$ which consists of the secants to $U$
through $(1,0,0)$ together with the line $[1,0,0]$. Each reported spread
$S'$ is described by the set of lines in $S' \setminus S$. Thus, $S'$
can be reconstructed by starting with the lines of $S$, discarding any
line of $S$ that meets one of the reported lines, then adding the
reported lines.

We represent the line $[x,y,z]$ in the form $(a,b,c)$ 
where $a$ (resp. $b$, $c$) is the discrete logarithm of $x$ 
(resp. $y$, $z$) with respect to $w$. * will be used to represent the
0 of the field. In addition to a representative spread, we also present a 
number of invariants of the spread, including:

\begin{itemize}
\item $|G|$ = order of the stabilizer of the spread within the automorphism 
group of the Hermitian unital.
\item Orbit = orbit structure of the stabilizer of the spread on its elements, 
reported as a set of pairs $(i,j)$ where $j$ is the number of orbits of size 
$i$. Pairs for which $j=0$ are not shown.
\item Type = set of pairs $(i,j)$ where $j$ is the number of points outside
of the unital $U$ that lie on exactly $i$ secants of the spread. Pairs for
which $j=0$ are not shown.
\item K = class of spread if generated from a known construction. Regular
spreads are denoted with an ``r'', cyclic spreads with a ``c'', the Andr\'e spreads
of Baker, et. al.~\cite{beks} with an ``a'', and non-Andr\'e subregular spreads, 
as defined in \cite{dover:subreg}, with an ``s''.
\end{itemize}

\small
\begin{center}
\begin{longtable}{|c|p{1.6in}|c|p{0.7in}|p{0.9in}|c|}
\caption[]{The spreads of the Hermitian unital in
$PG(2,16)$\label{PG216table}}\\
\hline
Idx & Spread & $|G|$ & Orbit & Type & K\\ \hline
\endfirsthead
\hline
Idx & Spread & $|G|$ & Orbit & Type & K\\ \hline
\endhead
\hline
\multicolumn{5}{r}{\emph {Continued on next page.}}
\endfoot
\hline
\endlastfoot
1 & (0,1,0), (0,13,12), (0,13,*), (0,1,*), (0,0,3), (0,6,8), (0,3,14), (0,2,8),
(0,12,6), (0,9,9), (0,2,14) & 156 & (13,1) & (0,130), (2,78) &  c \\ \hline
2 & (0,12,10), (0,6,4), (0,3,1), (0,0,13), (0,9,7) & 100 & (1,1), (2,1), (10,1)
& (0,100), (1,70), (2,36), (7,2) & a \\ \hline
3 &  & 1200 & (1,1), (12,1) & (0,75), (1,120), (2,12), (12,1) & r \\ \hline
\end{longtable}
\end{center}
\normalsize

To classify our spreads into the known types, we can use the following rules.
Regular spreads are distinguished by having a point off the unital lying on
$q^2-q$ lines of of the spread. From Baker, et. al.,~\cite{beks} cyclic
spreads exist only when $q$ is even, so we need only check spreads for even
$q$ with a transitive automorphism group. An algorithm to identify
subregular spreads is given in Dover~\cite{dover:subreg}, as is an
algorithm to distinguish within this class the spreads obtainable from
the triply-ruled $T_a$'s.

\section{$PG(2,25)$}

Let $w$ be a primitive element of GF(25) with minimal polynomial $x^2+4x+2$ over 
$GF(5)$. Our search found ten nonisomorphic spreads of the Hermitian unital, given in 
Table~\ref{PG225table}.

\small
\begin{center}
\begin{longtable}{|c|p{1.8in}|c|p{0.3in}|p{1.1in}|c|}
\caption[]{The spreads of the Hermitian unital in 
$PG(2,25)$\label{PG225table}}\\
\hline
Idx & Spread & $|G|$ & Orbit & Type & K\\ \hline
\endfirsthead
\hline
Idx & Spread & $|G|$ & Orbit & Type & K\\ \hline
\endhead
\hline
\multicolumn{5}{r}{\emph {Continued on next page.}}
\endfoot
\hline
\endlastfoot
1 &  & 1440 & (1,1), (20,1) & (0,144), (1,360), (2,20), (20,1) & r \\ \hline
2 & (0,4,0), (0,20,16), (0,12,8), (0,16,12), (0,8,4), (0,0,20) & 72 & (1,3), 
(6,3) & (0,216), (1,216), (2,91), (8,1), (14,1) & a \\ \hline
3 & (0,10,18), (0,2,10), (0,14,22), (0,6,14), (0,22,6), (0,18,2) & 72 & (1,3), 
(6,3) & (0,216), (1,216), (2,91), (8,1), (14,1) & a \\ \hline
4 & (0,18,5), (0,2,13), (0,14,1), (0,10,21), (0,6,17), (0,22,9) & 72 & (1,3), 
(6,3) & (0,216), (1,216), (2,91), (8,1), (14,1) & a \\ \hline
5 & (0,15,*), (0,3,*), (0,*,17), (0,*,1), (0,23,*), (0,*,21), (0,*,5), (0,*,13),
(0,19,*), (0,*,9), (0,7,*), (0,11,*) & 216 & (3,1), (18,1) & (0,252), (1,144), 
(2,126), (8,3) & a \\ \hline
6 & (0,16,1), (0,16,5), (0,8,17), (0,20,9), (0,12,1), (0,0,13), (0,20,5), 
(0,0,9), (0,4,17), (0,4,13), (0,8,21), (0,12,21) & 24 & (1,1), (2,2), (4,1), 
(12,1) & (0,240), (1,190), (2,68), (3,24), (7,2), (8,1) & s \\ \hline
7 & (0,5,10), (0,1,6), (0,6,17), (0,17,22), (0,10,21), (0,2,13), (0,21,2), 
(0,22,9), (0,14,1), (0,9,14), (0,18,5), (0,13,18) & 432 & (3,1), (18,1) & 
(0,216), (1,252), (2,18), (3,36), (8,3) & a \\ \hline
8 & (0,5,12), (0,13,14), (0,17,0), (0,17,18), (0,1,2), (0,21,4), (0,9,16), 
(0,9,10), (0,1,8), (0,21,22), (0,13,20), (0,5,6) & 144 & (1,1), (2,1), (6,1), 
(12,1) & (0,216), (1,252), (2,18), (3,36), (8,3) & a \\ \hline
9 & (0,23,3), (0,19,11), (0,10,11), (0,3,19), (0,*,21), (0,18,19), (0,22,11), 
(0,*,13), (0,*,5), (0,2,3), (0,14,3), (0,7,11), (0,6,19), (0,15,19), (0,11,3) & 
432 & (9,1), (12,1) & (0,216), (1,252), (2,36), (4,9), (5,12) &  \\ \hline
10 & (0,17,4), (0,15,3), (0,16,5), (0,*,4), (0,15,15), (0,17,16), (0,16,*), 
(0,4,*), (0,5,4), (0,3,3), (0,4,5), (0,3,15), (0,4,17), (0,*,16), (0,5,16), 
(0,16,17) & 336 & (21,1) & (0,224), (1,252), (3,28), (4,21) &  \\ \hline
\end{longtable}
\end{center}
\normalsize

Only two of the spreads found in $PG(2,25)$ do not belong to known families, namely spreads
9 and 10. Spread 10 admits a transitive automorphism group on the spread which
contains an index 2 copy of $PSL(2,7)$. Based on the interesting group
structure, this spread has been previously constructed by Durante and
Penttila~\cite{penttila}. While a very interesting spread, it seems 
unlikely that it will generalize to an infinite family of spreads for general $q$.

On the other hand, spread 9 seems capable of generalization. Spread 9
has the interesting property that its automorphism group has
the same order as spread 7, which is obtained from the triply-ruled construction of
Baker, et. al.~\cite{beks} via
\begin{displaymath}
\{[1,0,0],[0,1,0],[0,0,1]\} \cup V_3 \cup H_1 \cup D_2,
\end{displaymath}
using the notation of Equation~\ref{tripleruling}. Moreover this automorphism group is
significantly larger than the automorphism groups for other spreads obtainable from
this construction, primarily due to the presence of automorphisms that permute the
lines $\{[1,0,0],[0,1,0],[0,0,1]\}$, as opposed to leaving each of them fixed. It is
possible that there is some connection between these two spreads.

\section{$PG(2,49)$}

Let $w$ be a primitive element of GF(49) with minimal polynomial $x^2+6x+3$ over 
$GF(7)$. Our search found eighty-one nonisomorphic spreads of the Hermitian unital in
$PG(2,49)$. Suprisingly, all but two of these spreads are obtainable
from known constructions:1 is regular, 42 are Andr\'e, and 36 are subregular, but not
Andr\'e.  The two remaining spreads are listed as numbers 58 and 81 in
Table~\ref{PG249table}.

\small
\begin{center}
\begin{longtable}{|c|p{2.5in}|c|p{0.3in}|p{0.9in}|c|}
\caption[]{The spreads of the Hermitian unital in
$PG(2,49)$\label{PG249table}}\\
\hline
Idx & Spread & $|G|$ & Orbit & Type & K\\ \hline
\endfirsthead
\hline
Idx & Spread & $|G|$ & Orbit & Type & K\\ \hline
\endhead
\hline
\multicolumn{5}{r}{\emph {Continued on next page.}}
\endfoot
\hline
\endlastfoot
1 & (0,25,30), (0,12,1), (0,26,15), (0,1,6), (0,44,33), (0,14,3),
(0,6,43), (0,31,36), (0,20,9), (0,18,7), (0,36,25), (0,42,31),
(0,37,42), (0,7,12), (0,8,45), (0,2,39), (0,38,27), (0,30,19),
(0,0,37), (0,32,21), (0,43,0), (0,13,18), (0,19,24), (0,24,13) &
128 & (1,3), (8,5) & (0,832), (1,848), (2,360), (3,64), (10,1),
(18,2) & a \\ \hline
2 & (0,9,44), (0,3,38), (0,15,2), (0,10,29), (0,23,10),
(0,41,28), (0,45,32), (0,35,22), (0,4,23), (0,34,5), (0,16,35),
(0,29,16), (0,17,4), (0,27,14), (0,11,46), (0,40,11), (0,5,40),
(0,46,17), (0,21,8), (0,39,26), (0,47,34), (0,28,47), (0,33,20),
(0,22,41) & 128 & (1,3), (8,5) & (0,832), (1,848), (2,360),
(3,64), (10,1), (18,2) & a \\ \hline
3 & (0,36,7), (0,10,29), (0,0,19), (0,4,23), (0,34,5), (0,16,35),
(0,10,45), (0,24,43), (0,16,3), (0,46,33), (0,22,9), (0,34,21),
(0,40,11), (0,4,39), (0,46,17), (0,12,31), (0,40,27), (0,6,25),
(0,28,15), (0,42,13), (0,30,1), (0,28,47), (0,22,41), (0,18,37) &
128 & (1,3), (8,5) & (0,832), (1,848), (2,360), (3,64), (10,1),
(18,2) & a \\ \hline
4 & (0,0,19), (0,30,1), (0,46,33), (0,7,42), (0,4,39), (0,10,45),
(0,43,30), (0,19,6), (0,16,3), (0,12,31), (0,34,21), (0,25,12),
(0,28,15), (0,6,25), (0,40,27), (0,22,9), (0,37,24), (0,42,13),
(0,13,0), (0,31,18), (0,18,37), (0,24,43), (0,1,36), (0,36,7) &
128 & (1,3), (8,5) & (0,832), (1,848), (2,360), (3,64), (10,1),
(18,2) & a \\ \hline
5 & (0,46,33), (0,39,10), (0,10,45), (0,4,39), (0,11,46),
(0,16,3), (0,21,40), (0,34,21), (0,15,34), (0,17,4), (0,3,22),
(0,41,28), (0,35,22), (0,33,4), (0,27,46), (0,28,15), (0,47,34),
(0,40,27), (0,45,16), (0,29,16), (0,22,9), (0,9,28), (0,23,10),
(0,5,40) & 128 & (1,3), (8,5) & (0,832), (1,848), (2,360),
(3,64), (10,1), (18,2) & a \\ \hline
6 & (0,0,19), (0,30,1), (0,46,33), (0,39,10), (0,4,39),
(0,10,45), (0,16,3), (0,21,40), (0,12,31), (0,34,21), (0,15,34),
(0,3,22), (0,33,4), (0,27,46), (0,28,15), (0,6,25), (0,40,27),
(0,45,16), (0,22,9), (0,42,13), (0,9,28), (0,18,37), (0,24,43),
(0,36,7) & 128 & (1,3), (8,5) & (0,832), (1,848), (2,360),
(3,64), (10,1), (18,2) & a \\ \hline
7 & (0,26,20), (0,27,45), (0,20,38), (0,14,8), (0,15,33),
(0,21,39), (0,20,14), (0,32,2), (0,38,32), (0,39,9), (0,38,8),
(0,45,15), (0,2,20), (0,14,32), (0,8,26), (0,8,2), (0,26,44),
(0,32,26), (0,2,44), (0,9,27), (0,44,38), (0,3,21), (0,33,3),
(0,44,14) & 128 & (1,3), (8,5) & (0,832), (1,848), (2,360),
(3,64), (10,1), (18,2) & a \\ \hline
8 & (0,46,11), (0,21,34), (0,26,7), (0,16,29), (0,4,17),
(0,14,43), (0,20,1), (0,8,37), (0,27,40), (0,28,41), (0,10,23),
(0,33,46), (0,44,25), (0,15,28), (0,45,10), (0,38,19), (0,9,22),
(0,39,4), (0,32,13), (0,3,16), (0,22,35), (0,2,31), (0,34,47),
(0,40,5) & 128 & (1,3), (8,5) & (0,832), (1,848), (2,360),
(3,64), (10,1), (18,2) & a \\ \hline
9 & (0,38,30), (0,20,12), (0,32,24), (0,26,18), (0,4,36),
(0,3,35), (0,9,41), (0,46,30), (0,21,5), (0,39,23), (0,22,6),
(0,14,6), (0,44,36), (0,2,42), (0,45,29), (0,40,24), (0,16,0),
(0,33,17), (0,15,47), (0,10,42), (0,27,11), (0,8,0), (0,28,12),
(0,34,18) & 128 & (1,3), (8,5) & (0,832), (1,848), (2,360),
(3,64), (10,1), (18,2) & a \\ \hline
10 & (0,23,35), (0,29,41), (0,9,21), (0,21,33), (0,39,3),
(0,11,47), (0,47,11), (0,41,29), (0,17,5), (0,27,39), (0,5,41),
(0,45,9), (0,47,35), (0,15,27), (0,11,23), (0,29,17), (0,17,29),
(0,35,23), (0,23,11), (0,35,47), (0,3,15), (0,33,45), (0,5,17),
(0,41,5) & 128 & (1,3), (8,5) & (0,832), (1,848), (2,360),
(3,64), (10,1), (18,2) & a \\ \hline
11 & (0,34,10), (0,40,16), (0,4,28), (0,46,46), (0,36,12),
(0,40,40), (0,22,46), (0,42,18), (0,22,22), (0,18,42), (0,16,40),
(0,16,16), (0,28,4), (0,46,22), (0,6,30), (0,28,28), (0,24,0),
(0,12,36), (0,4,4), (0,10,34), (0,34,34), (0,0,24), (0,10,10),
(0,30,6) & 256 & (1,1), (2,1), (8,1), (16,2) & (0,768), (1,1040),
(2,168), (3,128), (10,1), (18,2) & a \\ \hline
12 & (0,27,32), (0,39,44), (0,25,30), (0,9,14), (0,26,15),
(0,1,6), (0,44,33), (0,14,3), (0,31,36), (0,20,9), (0,45,2),
(0,37,42), (0,7,12), (0,15,20), (0,8,45), (0,21,26), (0,2,39),
(0,38,27), (0,33,38), (0,3,8), (0,32,21), (0,43,0), (0,13,18),
(0,19,24) & 256 & (1,1), (2,1), (8,1), (16,2) & (0,768),
(1,1040), (2,168), (3,128), (10,1), (18,2) & a \\ \hline
13 & (0,29,40), (0,14,25), (0,11,22), (0,41,4), (0,17,28),
(0,20,15), (0,38,1), (0,23,34), (0,26,21), (0,8,19), (0,32,27),
(0,47,10), (0,38,33), (0,35,46), (0,8,3), (0,44,7), (0,5,16),
(0,20,31), (0,14,9), (0,32,43), (0,2,13), (0,2,45), (0,44,39),
(0,26,37) & 256 & (1,1), (2,1), (8,1), (16,2) & (0,768),
(1,1040), (2,168), (3,128), (10,1), (18,2) & a \\ \hline
14 & (0,4,46), (0,38,8), (0,32,2), (0,46,40), (0,10,4), (0,8,26),
(0,9,27), (0,26,44), (0,16,10), (0,3,21), (0,22,16), (0,20,38),
(0,21,39), (0,44,14), (0,14,32), (0,45,15), (0,34,28), (0,28,22),
(0,33,3), (0,39,9), (0,27,45), (0,15,33), (0,40,34), (0,2,20) &
256 & (1,1), (2,1), (8,1), (16,2) & (0,768), (1,1040), (2,168),
(3,128), (10,1), (18,2) & a \\ \hline
15 & (0,16,33), (0,22,23), (0,47,16), (0,46,47), (0,35,4),
(0,16,17), (0,40,41), (0,29,46), (0,28,29), (0,46,15), (0,34,3),
(0,23,40), (0,10,11), (0,17,34), (0,40,9), (0,28,45), (0,5,22),
(0,34,35), (0,22,39), (0,41,10), (0,4,5), (0,10,27), (0,4,21),
(0,11,28) & 128 & (1,3), (8,5) & (0,768), (1,1040), (2,168),
(3,128), (10,1), (18,2) & a \\ \hline
16 & (0,31,36), (0,19,24), (0,43,0), (0,14,3), (0,8,45),
(0,20,9), (0,32,21), (0,38,27), (0,2,39), (0,25,30), (0,7,12),
(0,1,6), (0,13,18), (0,37,42), (0,26,15), (0,44,33) & 128 &
(1,3), (8,5) & (0,768), (1,976), (2,296), (3,64), (10,2), (26,1)
& a \\ \hline
17 & (0,10,29), (0,41,28), (0,23,10), (0,35,22), (0,4,23),
(0,34,5), (0,16,35), (0,29,16), (0,17,4), (0,11,46), (0,40,11),
(0,5,40), (0,46,17), (0,47,34), (0,28,47), (0,22,41) & 256 &
(1,1), (2,1), (8,1), (16,2) & (0,768), (1,976), (2,296), (3,64),
(10,2), (26,1) & a \\ \hline
18 & (0,8,3), (0,39,2), (0,2,45), (0,45,8), (0,27,38), (0,14,9),
(0,44,39), (0,20,15), (0,33,44), (0,32,27), (0,21,32), (0,9,20),
(0,3,14), (0,26,21), (0,38,33), (0,15,26) & 256 & (1,1), (2,1),
(8,1), (16,2) & (0,768), (1,976), (2,296), (3,64), (10,2), (26,1)
& a \\ \hline
19 & (0,8,27), (0,26,45), (0,32,3), (0,23,10), (0,41,28),
(0,35,22), (0,14,33), (0,29,16), (0,17,4), (0,11,46), (0,2,21),
(0,5,40), (0,44,15), (0,20,39), (0,38,9), (0,47,34) & 128 &
(1,3), (8,5) & (0,768), (1,976), (2,296), (3,64), (10,2), (26,1)
& a \\ \hline
20 & (0,37,14), (0,39,32), (0,21,14), (0,7,32), (0,1,26),
(0,19,44), (0,9,2), (0,25,2), (0,31,8), (0,15,8), (0,45,38),
(0,27,20), (0,3,44), (0,43,20), (0,13,38), (0,33,26) & 128 &
(1,3), (8,5) & (0,768), (1,976), (2,296), (3,64), (10,2), (26,1)
& a \\ \hline
21 & (0,3,22), (0,10,45), (0,16,3), (0,46,33), (0,27,46),
(0,22,9), (0,9,28), (0,34,21), (0,4,39), (0,40,27), (0,45,16),
(0,28,15), (0,21,40), (0,15,34), (0,39,10), (0,33,4) & 256 &
(1,1), (2,1), (8,1), (16,2) & (0,768), (1,976), (2,296), (3,64),
(10,2), (26,1) & a \\ \hline
22 & (0,5,21), (0,23,39), (0,17,33), (0,47,15), (0,9,17),
(0,29,45), (0,35,3), (0,3,11), (0,33,41), (0,11,27), (0,21,29),
(0,45,5), (0,41,9), (0,15,23), (0,39,47), (0,27,35) & 256 &
(1,1), (2,1), (8,1), (16,2) & (0,768), (1,976), (2,296), (3,64),
(10,2), (26,1) & a \\ \hline
23 & (0,30,10), (0,24,4), (0,45,15), (0,27,45), (0,39,9),
(0,6,34), (0,15,33), (0,12,40), (0,9,27), (0,18,46), (0,33,3),
(0,0,28), (0,21,39), (0,42,22), (0,3,21), (0,36,16) & 32 & (1,1),
(2,7), (4,3), (16,1) & (0,768), (1,976), (2,296), (3,64), (10,2),
(26,1) & s \\ \hline
24 & (0,26,20), (0,27,45), (0,14,8), (0,15,33), (0,21,39),
(0,20,14), (0,38,32), (0,39,9), (0,45,15), (0,8,2), (0,32,26),
(0,2,44), (0,9,27), (0,3,21), (0,44,38), (0,33,3) & 128 & (1,3),
(8,5) & (0,768), (1,976), (2,296), (3,64), (10,2), (26,1) & a \\
\hline
25 & (0,21,34), (0,26,7), (0,14,43), (0,20,1), (0,8,37),
(0,27,40), (0,33,46), (0,44,25), (0,15,28), (0,45,10), (0,38,19),
(0,9,22), (0,39,4), (0,32,13), (0,3,16), (0,2,31) & 128 & (1,3),
(8,5) & (0,768), (1,976), (2,296), (3,64), (10,2), (26,1) & a \\
\hline
26 & (0,46,45), (0,10,41), (0,28,11), (0,34,33), (0,22,5),
(0,34,17), (0,40,39), (0,28,27), (0,40,23), (0,22,21), (0,4,35),
(0,10,9), (0,16,15), (0,46,29), (0,4,3), (0,16,47) & 128 & (1,3),
(8,5) & (0,768), (1,976), (2,296), (3,64), (10,2), (26,1) & a \\
\hline
27 & (0,1,6), (0,25,30), (0,13,18), (0,37,42), (0,19,24),
(0,43,0), (0,7,12), (0,31,36) & 128 & (1,3), (8,5) & (0,640),
(1,1168), (2,297), (10,1), (34,1) & a \\ \hline
28 & (0,35,22), (0,11,46), (0,5,40), (0,29,16), (0,17,4),
(0,23,10), (0,47,34), (0,41,28) & 128 & (1,3), (8,5) & (0,640),
(1,1168), (2,297), (10,1), (34,1) & a \\ \hline
29 & (0,28,15), (0,46,33), (0,16,3), (0,34,21), (0,10,45),
(0,4,39), (0,40,27), (0,22,9) & 128 & (1,3), (8,5) & (0,640),
(1,1168), (2,297), (10,1), (34,1) & a \\ \hline
30 & (0,28,16), (0,46,34), (0,16,4), (0,40,28), (0,10,46),
(0,34,22), (0,22,10), (0,4,40) & 128 & (1,3), (8,5) & (0,640),
(1,1168), (2,297), (10,1), (34,1) & a \\ \hline
31 & (0,45,15), (0,27,45), (0,39,9), (0,15,33), (0,9,27),
(0,33,3), (0,21,39), (0,3,21) & 128 & (1,3), (8,5) & (0,640),
(1,1168), (2,297), (10,1), (34,1) & a \\ \hline
32 & (0,27,32), (0,39,44), (0,25,30), (0,9,14), (0,1,6),
(0,31,36), (0,45,2), (0,37,42), (0,7,12), (0,15,20), (0,21,26),
(0,33,38), (0,3,8), (0,43,0), (0,13,18), (0,19,24) & 128 & (1,3),
(8,5) & (0,768), (1,912), (2,425), (18,1), (26,1) & a \\ \hline
33 & (0,9,44), (0,3,38), (0,15,2), (0,23,10), (0,41,28),
(0,45,32), (0,35,22), (0,29,16), (0,17,4), (0,27,14), (0,11,46),
(0,5,40), (0,21,8), (0,39,26), (0,47,34), (0,33,20) & 128 &
(1,3), (8,5) & (0,768), (1,912), (2,425), (18,1), (26,1) & a \\
\hline
34 & (0,23,10), (0,41,28), (0,35,22), (0,10,45), (0,29,16),
(0,17,4), (0,46,33), (0,16,3), (0,11,46), (0,22,9), (0,34,21),
(0,4,39), (0,5,40), (0,40,27), (0,28,15), (0,47,34) & 128 &
(1,3), (8,5) & (0,768), (1,912), (2,425), (18,1), (26,1) & a \\
\hline
35 & (0,9,44), (0,3,38), (0,15,2), (0,45,32), (0,10,45),
(0,16,3), (0,27,14), (0,46,33), (0,22,9), (0,34,21), (0,4,39),
(0,40,27), (0,21,8), (0,28,15), (0,39,26), (0,33,20) & 128 &
(1,3), (8,5) & (0,768), (1,912), (2,425), (18,1), (26,1) & a \\
\hline
36 & (0,17,45), (0,9,37), (0,47,27), (0,27,7), (0,41,21),
(0,15,43), (0,21,1), (0,35,15), (0,39,19), (0,33,13), (0,5,33),
(0,11,39), (0,45,25), (0,3,31), (0,29,9), (0,23,3) & 128 & (1,3),
(8,5) & (0,768), (1,912), (2,425), (18,1), (26,1) & a \\ \hline
37 & (0,6,24), (0,24,42), (0,36,6), (0,27,45), (0,15,33),
(0,21,39), (0,30,0), (0,39,9), (0,0,18), (0,45,15), (0,12,30),
(0,42,12), (0,18,36), (0,9,27), (0,3,21), (0,33,3) & 128 & (1,3),
(8,5) & (0,768), (1,912), (2,425), (18,1), (26,1) & a \\ \hline
38 & (0,29,23), (0,9,3), (0,41,35), (0,39,33), (0,17,11),
(0,5,47), (0,23,17), (0,21,15), (0,3,45), (0,11,5), (0,45,39),
(0,47,41), (0,33,27), (0,35,29), (0,27,21), (0,15,9) & 128 &
(1,3), (8,5) & (0,768), (1,912), (2,425), (18,1), (26,1) & a \\
\hline
39 & (0,20,38), (0,27,45), (0,15,33), (0,21,39), (0,32,2),
(0,39,9), (0,38,8), (0,45,15), (0,2,20), (0,14,32), (0,8,26),
(0,26,44), (0,9,27), (0,3,21), (0,33,3), (0,44,14) & 128 & (1,3),
(8,5) & (0,768), (1,912), (2,425), (18,1), (26,1) & a \\ \hline
40 & (0,33,14), (0,26,7), (0,9,38), (0,14,43), (0,20,1),
(0,8,37), (0,15,44), (0,44,25), (0,3,32), (0,38,19), (0,45,26),
(0,32,13), (0,2,31), (0,27,8), (0,39,20), (0,21,2) & 128 & (1,3),
(8,5) & (0,768), (1,912), (2,425), (18,1), (26,1) & a \\ \hline
41 & (0,7,*), (0,16,*), (0,1,*), (0,25,*), (0,19,*), (0,28,*),
(0,22,*), (0,31,*), (0,40,*), (0,13,*), (0,4,*), (0,46,*),
(0,34,*), (0,10,*), (0,37,*), (0,43,*) & 128 & (1,3), (8,5) &
(0,768), (1,912), (2,425), (18,1), (26,1) & a \\ \hline
42 & (0,31,36), (0,19,24), (0,13,30), (0,43,0), (0,37,6),
(0,43,12), (0,31,0), (0,7,24), (0,19,36), (0,1,18), (0,1,6),
(0,25,30), (0,7,12), (0,25,42), (0,13,18), (0,37,42) & 32 &
(1,1), (2,5), (4,4), (16,1) & (0,784), (1,928), (2,344), (3,48),
(10,2), (26,1) & s \\ \hline
43 & (0,16,4), (0,39,2), (0,4,40), (0,10,46), (0,33,44),
(0,22,10), (0,15,26), (0,9,20), (0,45,8), (0,46,34), (0,28,16),
(0,3,14), (0,21,32), (0,27,38), (0,34,22), (0,40,28) & 32 &
(1,1), (2,7), (4,3), (16,1) & (0,784), (1,928), (2,344), (3,48),
(10,2), (26,1) & s \\ \hline
44 & (0,41,5), (0,39,2), (0,45,8), (0,27,38), (0,11,23),
(0,33,44), (0,17,29), (0,21,32), (0,9,20), (0,23,35), (0,29,41),
(0,3,14), (0,35,47), (0,47,11), (0,15,26), (0,5,17) & 16 & (1,9),
(2,9), (8,2) & (0,784), (1,928), (2,344), (3,48), (10,2), (26,1)
& s \\ \hline
45 & (0,17,29), (0,41,5), (0,33,28), (0,5,17), (0,11,23),
(0,27,22), (0,21,16), (0,23,35), (0,15,10), (0,35,47), (0,47,11),
(0,3,46), (0,45,40), (0,9,4), (0,39,34), (0,29,41) & 32 & (1,1),
(2,7), (4,3), (16,1) & (0,784), (1,928), (2,344), (3,48), (10,2),
(26,1) & s \\ \hline
46 & (0,26,7), (0,27,45), (0,15,33), (0,14,43), (0,20,1),
(0,8,37), (0,21,39), (0,39,9), (0,45,15), (0,44,25), (0,38,19),
(0,32,13), (0,2,31), (0,9,27), (0,3,21), (0,33,3) & 32 & (1,3),
(2,4), (4,4), (16,1) & (0,784), (1,928), (2,344), (3,48), (10,2),
(26,1) & s \\ \hline
47 & (0,21,37), (0,22,16), (0,16,10), (0,34,28), (0,33,1),
(0,3,19), (0,28,22), (0,15,31), (0,46,40), (0,9,25), (0,4,46),
(0,27,43), (0,45,13), (0,10,4), (0,39,7), (0,40,34) & 16 & (1,5),
(2,11), (8,2) & (0,784), (1,928), (2,344), (3,48), (10,2), (26,1)
& s \\ \hline
48 & (0,26,6), (0,44,24), (0,31,39), (0,37,45), (0,43,3),
(0,32,12), (0,19,27), (0,13,21), (0,1,9), (0,38,18), (0,25,33),
(0,2,30), (0,8,36), (0,20,0), (0,14,42), (0,7,15) & 16 & (1,5),
(2,11), (8,2) & (0,784), (1,928), (2,344), (3,48), (10,2), (26,1)
& s \\ \hline
49 & (0,7,42), (0,11,46), (0,43,30), (0,8,27), (0,19,6),
(0,20,39), (0,26,45), (0,17,4), (0,41,28), (0,35,22), (0,14,33),
(0,25,12), (0,38,9), (0,47,34), (0,29,16), (0,37,24), (0,13,0),
(0,31,18), (0,23,10), (0,1,36), (0,32,3), (0,2,21), (0,5,40),
(0,44,15) & 128 & (1,3), (8,5) & (0,896), (1,656), (2,552),
(10,1), (18,2) & a \\ \hline
50 & (0,1,30), (0,37,2), (0,27,40), (0,43,8), (0,7,36),
(0,13,42), (0,25,6), (0,31,44), (0,37,18), (0,21,34), (0,9,22),
(0,19,0), (0,13,26), (0,15,28), (0,43,24), (0,3,16), (0,25,38),
(0,7,20), (0,1,14), (0,39,4), (0,31,12), (0,33,46), (0,45,10),
(0,19,32) & 128 & (1,3), (8,5) & (0,896), (1,656), (2,552),
(10,1), (18,2) & a \\ \hline
51 & (0,15,19), (0,3,7), (0,33,37), (0,41,28), (0,23,10),
(0,35,22), (0,29,16), (0,17,4), (0,11,46), (0,39,43), (0,45,1),
(0,5,40), (0,21,25), (0,47,34), (0,9,13), (0,27,31) & 32 & (1,1),
(2,5), (4,4), (16,1) & (0,800), (1,894), (2,362), (3,48), (9,2),
(26,1) & s \\ \hline
52 & (0,36,7), (0,0,19), (0,10,45), (0,24,43), (0,16,3),
(0,46,33), (0,22,9), (0,34,21), (0,4,39), (0,12,31), (0,40,27),
(0,6,25), (0,28,15), (0,42,13), (0,30,1), (0,18,37) & 128 &
(1,3), (8,5) & (0,832), (1,784), (2,488), (10,2), (26,1) & a \\
\hline
53 & (0,46,11), (0,26,7), (0,16,29), (0,4,17), (0,14,43),
(0,20,1), (0,8,37), (0,28,41), (0,10,23), (0,44,25), (0,38,19),
(0,32,13), (0,22,35), (0,2,31), (0,34,47), (0,40,5) & 128 &
(1,3), (8,5) & (0,832), (1,784), (2,488), (10,2), (26,1) & a \\
\hline
54 & (0,8,3), (0,2,45), (0,16,*), (0,3,47), (0,14,9), (0,20,15),
(0,44,39), (0,28,*), (0,22,*), (0,9,5), (0,15,11), (0,45,41),
(0,32,27), (0,40,*), (0,21,17), (0,4,*), (0,46,*), (0,34,*),
(0,27,23), (0,33,29), (0,38,33), (0,26,21), (0,10,*), (0,39,35) &
16 & (1,5), (2,7), (8,1), (16,1) & (0,864), (1,814), (2,345),
(3,64), (4,16), (9,2), (10,1), (18,1) & s \\ \hline
55 & (0,39,2), (0,33,44), (0,20,15), (0,20,1), (0,26,7),
(0,15,26), (0,26,21), (0,32,27), (0,9,20), (0,45,8), (0,38,33),
(0,8,3), (0,14,9), (0,3,14), (0,14,43), (0,38,19), (0,21,32),
(0,44,25), (0,32,13), (0,27,38), (0,8,37), (0,2,45), (0,2,31),
(0,44,39) & 32 & (1,3), (2,4), (4,2), (8,1), (16,1) & (0,864),
(1,814), (2,345), (3,64), (4,16), (9,2), (10,1), (18,1) & s \\
\hline
56 & (0,9,4), (0,39,34), (0,21,16), (0,27,22), (0,23,35),
(0,3,46), (0,15,10), (0,29,41), (0,30,46), (0,47,11), (0,42,10),
(0,0,16), (0,36,4), (0,6,22), (0,33,28), (0,12,28), (0,11,23),
(0,45,40), (0,18,34), (0,17,29), (0,35,47), (0,5,17), (0,24,40),
(0,41,5) & 16 & (1,7), (2,6), (8,1), (16,1) & (0,864), (1,814),
(2,345), (3,64), (4,16), (9,2), (10,1), (18,1) & s \\ \hline
57 & (0,8,3), (0,2,45), (0,3,47), (0,14,9), (0,20,15), (0,44,39),
(0,9,5), (0,15,11), (0,45,41), (0,32,27), (0,21,17), (0,27,23),
(0,33,29), (0,38,33), (0,26,21), (0,39,35) & 16 & (1,5), (2,11),
(8,2) & (0,784), (1,935), (2,329), (3,56), (9,1), (10,1), (26,1)
& s \\ \hline
58 & (0,6,20), (0,3,4), (0,20,42), (0,32,6), (0,15,16),
(0,28,45), (0,40,9), (0,16,3), (0,9,28), (0,4,39), (0,44,18),
(0,40,27), (0,28,15), (0,45,16), (0,42,8), (0,27,28), (0,30,44),
(0,21,40), (0,4,21), (0,39,40), (0,16,33), (0,18,32), (0,33,4),
(0,8,30) & 192 & (1,1), (6,1), (12,1), (24,1) & (0,864), (1,792),
(2,420), (5,24), (6,6), (18,1) & \\ \hline
59 &  & 5376 & (1,1), (42,1) & (0,384), (1,1680), (2,42), (42,1)
& r \\ \hline
60 & (0,21,34), (0,17,45), (0,47,27), (0,27,40), (0,41,21),
(0,33,46), (0,35,15), (0,45,10), (0,15,28), (0,39,4), (0,5,33),
(0,9,22), (0,11,39), (0,3,16), (0,29,9), (0,23,3) & 16 & (1,5),
(2,11), (8,2) & (0,792), (1,911), (2,353), (3,48), (9,1), (10,1),
(26,1) & s \\ \hline
61 & (0,46,9), (0,26,7), (0,34,45), (0,40,3), (0,14,43),
(0,10,21), (0,20,1), (0,8,37), (0,16,27), (0,44,25), (0,38,19),
(0,28,39), (0,32,13), (0,4,15), (0,2,31), (0,22,33) & 16 & (1,7),
(2,10), (8,2) & (0,792), (1,911), (2,353), (3,48), (9,1), (10,1),
(26,1) & s \\ \hline
62 & (0,39,2), (0,26,7), (0,45,8), (0,27,38), (0,14,43),
(0,20,1), (0,8,37), (0,33,44), (0,21,32), (0,9,20), (0,44,25),
(0,38,19), (0,32,13), (0,2,31), (0,3,14), (0,15,26) & 16 & (1,7),
(2,10), (8,2) & (0,792), (1,911), (2,353), (3,48), (9,1), (10,1),
(26,1) & s \\ \hline
63 & (0,9,41), (0,46,11), (0,27,11), (0,33,17), (0,16,29),
(0,4,17), (0,39,23), (0,28,41), (0,10,23), (0,15,47), (0,21,5),
(0,22,35), (0,34,47), (0,3,35), (0,45,29), (0,40,5) & 16 & (1,5),
(2,11), (8,2) & (0,792), (1,911), (2,353), (3,48), (9,1), (10,1),
(26,1) & s \\ \hline
64 & (0,4,15), (0,39,2), (0,22,33), (0,16,27), (0,33,44),
(0,26,7), (0,20,1), (0,15,26), (0,40,3), (0,9,20), (0,45,8),
(0,3,14), (0,34,45), (0,14,43), (0,38,19), (0,21,32), (0,44,25),
(0,10,21), (0,32,13), (0,27,38), (0,28,39), (0,8,37), (0,46,9),
(0,2,31) & 16 & (1,7), (2,6), (8,3) & (0,848), (1,815), (2,377),
(3,64), (10,1), (17,1), (18,1) & s \\ \hline
65 & (0,9,4), (0,4,15), (0,39,34), (0,21,16), (0,27,22),
(0,3,46), (0,15,10), (0,22,33), (0,16,27), (0,26,7), (0,20,1),
(0,40,3), (0,33,28), (0,45,40), (0,34,45), (0,14,43), (0,38,19),
(0,10,21), (0,44,25), (0,32,13), (0,28,39), (0,8,37), (0,46,9),
(0,2,31) & 32 & (1,3), (2,4), (4,2), (8,1), (16,1) & (0,864),
(1,822), (2,321), (3,88), (4,8), (9,2), (10,1), (18,1) & s \\
\hline
66 & (0,47,15), (0,39,8), (0,15,32), (0,35,3), (0,16,41),
(0,23,39), (0,40,17), (0,10,35), (0,29,45), (0,17,33), (0,11,27),
(0,41,9), (0,34,11), (0,28,5), (0,5,21), (0,4,29), (0,46,23),
(0,9,26), (0,45,14), (0,22,47), (0,3,20), (0,21,38), (0,33,2),
(0,27,44) & 8 & (1,19), (8,3) & (0,864), (1,822), (2,321),
(3,88), (4,8), (9,2), (10,1), (18,1) & s \\ \hline
67 & (0,8,36), (0,30,12), (0,24,6), (0,14,42), (0,33,47),
(0,0,30), (0,2,30), (0,18,0), (0,27,41), (0,6,36), (0,32,12),
(0,26,6), (0,39,5), (0,36,18), (0,12,42), (0,42,24), (0,20,0),
(0,45,11), (0,3,17), (0,9,23), (0,38,18), (0,15,29), (0,44,24),
(0,21,35) & 8 & (1,19), (8,3) & (0,864), (1,822), (2,321),
(3,88), (4,8), (9,2), (10,1), (18,1) & s \\ \hline
68 & (0,42,22), (0,26,7), (0,36,16), (0,18,46), (0,27,45),
(0,15,33), (0,14,43), (0,20,1), (0,8,37), (0,21,39), (0,30,10),
(0,6,34), (0,24,4), (0,39,9), (0,45,15), (0,44,25), (0,38,19),
(0,32,13), (0,2,31), (0,12,40), (0,9,27), (0,3,21), (0,33,3),
(0,0,28) & 16 & (1,7), (2,6), (8,1), (16,1) & (0,848), (1,848),
(2,327), (3,64), (4,16), (10,3), (18,1) & s \\ \hline
69 & (0,37,2), (0,6,34), (0,42,22), (0,43,8), (0,12,40),
(0,9,27), (0,30,10), (0,31,44), (0,3,21), (0,24,4), (0,13,26),
(0,21,39), (0,45,15), (0,36,16), (0,25,38), (0,18,46), (0,33,3),
(0,39,9), (0,0,28), (0,27,45), (0,7,20), (0,1,14), (0,15,33),
(0,19,32) & 48 & (1,1), (3,6), (24,1) & (0,816), (1,960),
(2,183), (3,144), (10,3), (18,1) & s \\ \hline
70 & (0,27,33), (0,24,16), (0,22,42), (0,36,28), (0,3,9),
(0,15,21), (0,33,39), (0,16,36), (0,9,15), (0,10,30), (0,30,22),
(0,21,27), (0,46,18), (0,0,40), (0,39,45), (0,28,0), (0,12,4),
(0,4,24), (0,18,10), (0,42,34), (0,40,12), (0,6,46), (0,45,3),
(0,34,6) & 96 & (1,1), (3,2), (6,2), (24,1) & (0,816), (1,960),
(2,183), (3,144), (10,3), (18,1) & s \\ \hline
71 & (0,46,10), (0,40,4), (0,30,46), (0,16,28), (0,22,34),
(0,20,15), (0,42,10), (0,0,16), (0,34,46), (0,36,4), (0,26,21),
(0,32,27), (0,6,22), (0,38,33), (0,12,28), (0,8,3), (0,18,34),
(0,14,9), (0,10,22), (0,24,40), (0,4,16), (0,2,45), (0,28,40),
(0,44,39) & 16 & (1,7), (2,6), (8,1), (16,1) & (0,856), (1,846),
(2,297), (3,96), (4,8), (9,2), (10,1), (18,1) & s \\ \hline
72 & (0,1,33), (0,8,36), (0,31,15), (0,11,40), (0,14,42),
(0,2,30), (0,5,34), (0,13,45), (0,43,27), (0,47,28), (0,23,4),
(0,32,12), (0,41,22), (0,26,6), (0,29,10), (0,17,46), (0,7,39),
(0,19,3), (0,35,16), (0,20,0), (0,25,9), (0,38,18), (0,37,21),
(0,44,24) & 16 & (1,5), (2,7), (8,1), (16,1) & (0,856), (1,846),
(2,297), (3,96), (4,8), (9,2), (10,1), (18,1) & s \\ \hline
73 & (0,4,15), (0,46,10), (0,22,33), (0,16,27), (0,40,4),
(0,30,46), (0,16,28), (0,22,34), (0,42,10), (0,0,16), (0,34,46),
(0,36,4), (0,6,22), (0,40,3), (0,12,28), (0,18,34), (0,34,45),
(0,10,22), (0,10,21), (0,24,40), (0,28,39), (0,4,16), (0,46,9),
(0,28,40) & 8 & (1,19), (8,3) & (0,872), (1,798), (2,345),
(3,80), (4,8), (9,2), (10,1), (18,1) & s \\ \hline
74 & (0,36,1), (0,8,36), (0,5,45), (0,0,13), (0,11,3), (0,14,42),
(0,23,15), (0,35,27), (0,41,33), (0,2,30), (0,12,25), (0,6,19),
(0,32,12), (0,29,21), (0,26,6), (0,30,43), (0,17,9), (0,18,31),
(0,47,39), (0,20,0), (0,24,37), (0,38,18), (0,42,7), (0,44,24) &
16 & (1,5), (2,7), (8,1), (16,1) & (0,872), (1,798), (2,345),
(3,80), (4,8), (9,2), (10,1), (18,1) & s \\ \hline
75 & (0,18,34), (0,42,22), (0,6,22), (0,36,16), (0,18,46),
(0,24,40), (0,27,45), (0,42,10), (0,30,46), (0,15,33), (0,21,39),
(0,36,4), (0,30,10), (0,6,34), (0,12,28), (0,24,4), (0,39,9),
(0,45,15), (0,0,16), (0,12,40), (0,9,27), (0,3,21), (0,33,3),
(0,0,28) & 32 & (1,3), (2,4), (4,2), (8,1), (16,1) & (0,912),
(1,672), (2,471), (3,48), (10,3), (18,1) & s \\ \hline
76 & (0,4,46), (0,16,4), (0,39,2), (0,46,40), (0,4,40), (0,10,4),
(0,10,46), (0,33,44), (0,22,10), (0,15,26), (0,9,20), (0,16,10),
(0,22,16), (0,45,8), (0,46,34), (0,28,16), (0,34,28), (0,3,14),
(0,28,22), (0,21,32), (0,27,38), (0,34,22), (0,40,34), (0,40,28)
& 48 & (1,1), (3,4), (6,1), (24,1) & (0,864), (1,816), (2,327),
(3,96), (10,3), (18,1) & s \\ \hline
77 & (0,9,4), (0,8,36), (0,39,34), (0,21,16), (0,27,22),
(0,23,35), (0,3,46), (0,15,10), (0,29,41), (0,14,42), (0,2,30),
(0,47,11), (0,32,12), (0,26,6), (0,33,28), (0,11,23), (0,45,40),
(0,17,29), (0,20,0), (0,35,47), (0,5,17), (0,38,18), (0,41,5),
(0,44,24) & 48 & (1,1), (3,4), (6,1), (24,1) & (0,864), (1,816),
(2,327), (3,96), (10,3), (18,1) & s \\ \hline
78 & (0,4,46), (0,39,2), (0,23,35), (0,29,41), (0,46,40),
(0,10,4), (0,33,44), (0,47,11), (0,15,26), (0,9,20), (0,16,10),
(0,22,16), (0,45,8), (0,11,23), (0,34,28), (0,3,14), (0,17,29),
(0,21,32), (0,28,22), (0,35,47), (0,27,38), (0,5,17), (0,40,34),
(0,41,5) & 16 & (1,5), (2,7), (8,1), (16,1) & (0,856), (1,832),
(2,327), (3,80), (4,8), (10,3), (18,1) & s \\ \hline
79 & (0,4,46), (0,9,4), (0,39,34), (0,21,16), (0,27,22),
(0,23,35), (0,3,46), (0,15,10), (0,29,41), (0,46,40), (0,10,4),
(0,47,11), (0,33,28), (0,16,10), (0,22,16), (0,11,23), (0,45,40),
(0,34,28), (0,17,29), (0,35,47), (0,28,22), (0,5,17), (0,40,34),
(0,41,5) & 16 & (1,5), (2,7), (8,1), (16,1) & (0,864), (1,800),
(2,375), (3,48), (4,16), (10,3), (18,1) & s \\ \hline
80 & (0,11,31), (0,1,33), (0,31,15), (0,23,43), (0,29,1),
(0,11,40), (0,47,19), (0,35,7), (0,5,34), (0,41,13), (0,13,45),
(0,43,27), (0,47,28), (0,23,4), (0,41,22), (0,29,10), (0,17,37),
(0,17,46), (0,7,39), (0,19,3), (0,35,16), (0,25,9), (0,5,25),
(0,37,21) & 16 & (1,5), (2,7), (8,1), (16,1) & (0,880), (1,782),
(2,345), (3,96), (9,2), (10,1), (18,1) & s \\ \hline
81 & (0,38,32), (0,2,44), (0,3,31), (0,11,45), (0,45,5),
(0,9,17), (0,23,3), (0,16,36), (0,26,20), (0,47,33), (0,15,43),
(0,25,12), (0,35,15), (0,21,29), (0,23,9), (0,28,0), (0,4,24),
(0,39,19), (0,37,24), (0,35,21), (0,13,0), (0,40,12), (0,1,36),
(0,47,27), (0,11,39), (0,27,7), (0,33,41), (0,14,8) & 32 & (1,3),
(4,6), (16,1) & (0,896), (1,760), (2,372), (3,48), (4,16), (5,8),
(6,5), (10,1), (14,1) &  \\ \hline
\end{longtable}
\end{center}
\normalsize

\section{Conclusion}
There are several potentially interesting directions for further research based upon
the results here. As mentioned previously,
it is likely possible to port the existing search algorithms to more robust hardware
 and complete a similar search for spreads of the Hermitian
unital in $PG(2,64)$. Such a search would likely provide new examples of spreads 
for even $q$.

In addition to spread problems, the problem of finding resolutions of Hermitian unitals,
namely partitions of the blocks of a unital into pairwise disjoint spreads, is of some
geometric interest. As one application, the author~\cite{thesis} provides a construction
of a fan of the three-dimensional Hermitian variety from a resolution of the Hermitian
unital. Currently, only the obvious resolution of the Hermitian unital into regular 
spreads whose starting points are all collinear is known.


\begin{thebibliography}{1}

\bibitem{beks} {\bf R.D.~Baker, G.L.~Ebert, G.~Korchm\'aros, T.~Sz\"onyi}:
Orthogonally divergent spreads of Hermitian curves,
In {\em Finite Geometry and Combinatorics}, volume 191 of {\em London Math. Soc. 
Lecture Note Ser.}, pages 17--30, Cambridge University Press, Cambridge, 1993 (Deinze, 1992).

\bibitem{cannon} {\bf J.~Cannon and C.~Playoust}:
{\em An Introduction to MAGMA},
University of Sydney, Sydney, 1993.

\bibitem{dover:ree} {\bf Jeremy M. Dover}:
Spreads and resolutions of Ree unitals,
{\em Ars Combin.} {\bf 54} (2000), \mbox(301--309).

\bibitem{dover:subreg} {\bf Jeremy M. Dover}:
Subregular spreads of Hermitian unitals,
{\em Des. Codes Crypt.} {\bf 39} (2006).

\bibitem{thesis} {\bf Jeremy M. Dover}:
{\em Theory and Applications of Spreads of Geometric Spaces},
Ph.D. Thesis, University of Delaware, 1996.

\bibitem{penttila}{\bf  N.~Durante and T.~Penttila}:
\newblock private communication.

\end{thebibliography}
\end{document}